\documentclass[12pt]{article}
\usepackage{amsthm}
\usepackage{latexsym}
\usepackage{fleqn}
\usepackage[latin1]{inputenc}
\usepackage{amssymb}
\usepackage{amsfonts}

\newcommand{\bR}{{\mathbb{R}}}
\newcommand{\bN}{{\mathbb{N}}}

\newcommand{\cC}{{\mathcal C}}

\usepackage[T1]{fontenc}           
\usepackage{textcomp}              
\usepackage{epsfig}        

\newtheorem{theorem}{Theorem}[section]
\newtheorem{lemma}{Lemma}[section]
\newtheorem{corollary}{Corollary}[section]

\theoremstyle{definition} 
\newtheorem{definition}{Definition}[section]  

\newtheorem{remark}{Remark}[section]
\newtheorem{example}{Example}[section]

\begin{document}

\begin{center}
\textbf{\Large On the polygonal diameter of the interior, resp. exterior, of a simple closed polygon in the plane}\\[0.5cm]

\emph{Yaakov S. Kupitz\footnote{Partially supported by the Landau
Center at the Mathematics Institute of the Hebrew University of
Jerusalem $($supported by Minerva Foundation, Germany$)$, and by Deutsche
Forschungsgemeinschaft.}$\,^*$,  Horst Martini~$^{**}$,  Micha A. Perles~$^*$}

$^*$  Institute of Mathematics, The Hebrew University of Jerusalem, Jerusalem, ISRAEL, kupitz@math.huji.ac.il; perles@math.huji.ac.il\\[0.2cm]
$^{**}$ Faculty of Mathematics, University of Technology, 09107 Chemnitz, GERMANY, martini@mathematik.tu-chemnitz.de
\end{center}

\begin{abstract}

We give a tight upper bound on the polygonal diameter of the interior, resp. exterior, of a simple $n$-gon, $n \ge 3$,
in the plane as a function of $n$, and describe an $n$-gon $(n \ge 3)$ for which both upper bounds (for the interior and the exterior)
are attained \emph{simultaneously}.
\end{abstract}



\section{Introduction}

The following is well known

\begin{theorem}{\rm (The Jordan theorem)}\label{theo1.1}
Let $f:[0,1] \to \bR^2$ be a simple closed curve in the plane ($f$ is continous, $f(0) = f(1)$ and $f(u) \not= f(v)$ for $0 < u < v \le 1$).
Define $P=_{\rm def}$ {\rm image}$f= \{f(u) : 0 \le u \le 1\}$, the image of $f$. Then $\bR^2 \setminus P = U_0 \cup U_1$, where $U_0, U_1$ are connected open,
non-empty mutually disjoint sets, $U_0$ is bounded (interior), $U_1$ is unbounded (exterior), and $P = {\rm bd}(U_0) = {\rm bd} (U_1)$.
\end{theorem}

The proof of this theorem is not easy; see \cite{Bertoglio}, \cite{Lawson}, \cite{Thomassen}, \cite[p. 37 ff.]{Moise}, \cite[vol. I, pp. 39-64]{Aleksandrov},
\cite[pp. 285 ff.]{Kuratowski}, and the survey \cite{Dostal}.
When the curve $P$ is polygonal, however, i.e., when $f$ is piecewise affine, the theorem becomes elementary:

\begin{theorem} {\rm(The piecewise affine Jordan theorem)} \label{theo1.2}
Let $p_0,p_1,\dots,p_{n-1},p_n=p_0, n \ge 3$, be ($n$ distinct) points in $\bR^2$. Assume that the polygon $P=_{\rm def} \bigcup\limits^n_{i=1}
[p_{i-1},p_i]$ is simple, i.e., the segments $[p_{i-1},p_i]$ do not intersect except for common endpoints: $\{p_i\} = [p_{i-1},p_i] \cap [p_i,p_{i+1}]$
for $1 \le i \le n-1, \{p_0\} = [p_0,p_1] \cap [p_{n-1},p_0]$. Then $\bR^2 \setminus P = U_0 \cup U_1$ with the same properties of $U_0,U_1$ listed
above {\rm (Theorem \ref{theo1.1})}.
\end{theorem}

\begin{definition}\label{def1.1}
A polygon $P$ satisfying the conditions of Theorem \ref{theo1.2} is a \emph{simple closed $n$-gon}. The bounded [resp. unbounded] domain $U_0$
[resp. $U_1$] is the \emph{interior} [resp. \emph{exterior}], denoted by int$P$ [resp. ext$P$], of $P$.
\end{definition}

A particularly simple proof of Theorem \ref{theo1.2} is known as  the ``raindrop proof'', see \cite[pp. 267-269]{Courant}, \cite[pp. 281-285]{Hille},
\cite[pp. 27-29]{Bensen}, or \cite[pp. 16-18]{Moise}. We reproduce this proof in a somewhat more complete and formal form than usually given in the literature for later
reference to some of its parts.

So we first prove Theorem \ref{theo1.2} (in Paragraphs 2 and 3 below). Then, squeezing this proof, a \emph{tight} upper bound on the polygonal diameter
of int$P$ [resp. ext$P$] (see Definition \ref{def3.2} below) is given as a function of $n$, and an $n$-gon $(n \ge 3)$ for which both upper bounds are
attained \emph{simultaneously} is described (see Theorem \ref{theo4.1} below). The $d$-dimensional analogue $(d \ge 2)$  of this problem was discussed
in \cite[Theorem 3.2]{Perles}. There we gave upper bounds on the polygonal diameter of int$\cC$, resp. ext$\cC$, for a polyhedral $(d-1)$-pseudomanifold $\cC$
in $\bR^d$ as a function of the number $n$ of its facets and $d$. The bounds given there are shown to be \emph{almost} tight (see \cite[Section 4]{Perles}), whereas the bounds given here (for $d = 2$) are tight. Another novelty of the present paper is that there is an $n$-gon $P$ in $\bR^2$ for which \emph{both} upper bounds (on the polygonal diameter of int$P$ and ext$P$) are attained
(simultanously), as said above, whereas for $d \ge 3$ the examples given in \cite[Section 4]{Perles} (namely one for int$\cC$ and another
one for ext$\cC$) are \emph{different} from each other.

For the sake of the proof of Theorem \ref{theo1.2}, we split it into two statements: Let $P$ be a simple closed polygon in $\bR^2$.

(E) (separation): $\bR^2 \setminus P$ is the disjoint union of two open sets, int$P$ and ext$P$. The boundary of each one of these sets
is $P$; int$P$ is bounded and ext$P$ is unbounded.

(F) (connectivity): The sets int$P$ and ext$P$ are [polygonally] connected.

We shall prove (E) (Paragraph 2) by constructing a continuous function $f: \bR^2 \setminus P \to \{0,1\}$ which attains both values $0$ and $1$ in every
neighborhood of every point $x \in P$, and defining ext$P = f^{-1}(0)$, int$P = f^{-1}(1)$. Statement (F) (polygonal connectivity of int$P$ and of ext$P$) follows
from Theorem \ref{theo3.1} below.

\section{A ``raindrop'' proof of (E)}

The construction of $f$ will be performed in three steps:

\textbf{Preliminary step:} Choosing a ``generic'' direction.

Choose an orthogonal basis $(u,v)$ for $\bR^2$ so that no two
vertices of $P$ have the same $x$-coordinate. Intuitively: the polygon $P$ is drawn as a paper; rotate the paper so that no two vertices
lie one above the other. Formally: let $L_1, \dots, L_t$ be all lines spanned by subsets of $\{p_1,\dots,p_n\}$. For $i=1,\dots,t$ let
$L^0_i =_{\rm def} L_i-L_i$ be the linear ($1$-dimensional) subspace parallel to $L_i$. Choose a unit vector $v \in \bR^2 \setminus \bigcup \limits^t_{i=1} L^0_i$
(``$v$'' for ``vertical'').
The vector $v$ is our direction ``up'', and $-v$ is pointing ``down''. By our choice of $v$, a line $L$, spanned by the vertices of $P$, will meet a line parallel
to $v$ in at most one point.

For a point $p \in \bR^2 \setminus P$ denote by $R(p)$ the closed vertical ``pointing down'' half-line $R(p) =_{\rm def}
\{p-\lambda v: 0 \le \lambda < \infty\}$. $R(p)$ is the path of a ``raindrop'' emanating from $p$. We divide $\bR^2 \setminus P$ into two disjoint sets
\[
\begin{array}{lll}
S_0 & =_{\rm def} & \{p \in \bR^2 \setminus P: R(p) \, \mbox{ does not meet any vertex of }\, P\}\,,\\
S_1 & =_{\rm def} & \{ p \in \bR^2 \setminus P: R(p) \, \mbox{ meets exactly one vertex of } \, P\}\,.
\end{array}
\]
(By our choice of $v$, we have $\bR^2 \setminus P = S_0 \cup S_1$.)
We shall define $f$ on $S_0$ (= Step I), then extend it (continuously) to $S_1$ (= Step II).
The following notation will be used: For a set $A \subset \bR^2$,
$A^+ =_{\rm def} \{a + \lambda v: a \in A, \lambda \ge 0\}$.

Thus $A^+$ is the set of points that lie
``above'' $A$. If $A$ is closed, then $A^+$ is closed. Note that (for all $p \in \bR^2$ and $A \subset \bR^2$):
\begin{equation}\label{eq1}
R(p) \, \mbox{ meets } \, A \, \mbox{ iff } \, p \in A^+\,.
\end{equation}

\textbf{Step I:} Define $f$ on $S_0$.

For $p \in S_o$ denote by $r(p)$ the number of edges of $P$ met by $R(p)$, and define $f(p) =_{\rm def} {\rm par} (r(p)) =_{\rm def} \frac{1}{2} (1-(-1)^{r(p)})$, the parity
of $r(p)$ ($f(p) = 0$ if $r(p)$ is even, $1$ if $r(p)$ is odd).

\begin{center}
\end{center}

\begin{center}
Fig. 1: the function $r(p)$ \hfill Fig. 2: the parity function $f(p) = {\rm par}(r(p))$
\end{center}

Next we show that $S_0$ is a dense open subset of $\bR^2$, and that $f: S_0 \to \{0,1\}$ is a continuous, hence locally constant function. Using vert$P$ for the set of vertices of $P$,
we have in view of (\ref{eq1})
\begin{equation}\label{eq2}
S_0 = \bR^2 \setminus (P \cup (\mbox{vert}P)^+)\,.
\end{equation}
The set $({\rm vert}P)^+$ is closed, same as $P$. Thus $S_0$ is an open subset of $\bR^2$. Moreover, the set $P \cup ({\rm vert}P)^+$ can be covered by a finite
number of lines in $\bR^2$. It follows that $S_0$ is dense in $R^2$.

Continuity of $f$: Assume $x \in S_0$. Let $\varepsilon$ be the (positive) distance from $x$ to $P \cup ({\rm vert}P)^+ (= \bR^2 \setminus S_0)$.
If $x' \in \bR^2, \|x-x'\| < \varepsilon$, then the segment $[x,x']$ does not meet $P \cup ({\rm vert}P)^+$. Let $e = [p_{i-1},p_i] \, (1 \le i \le n)$ be any
edge of $P$. The set $e^+$ is a closed, convex, unbounded and full-dimensional polyhedral subset of $\bR^2$, whose boundary consists of the lower edge $e$ and the side
edges $p^+_{i-1}, p^+_i$. Thus bd$e^+ \subset P \cup ({\rm vert} P)^+$, and therefore the segment $[x,x']$ does not meet the boundary of $e^+$. It follows that
$x' \in e^+$ iff $x \in e^+$, i.e., $R(x)$ meets $e$ iff $R(x')$ meets $e$. This is true for all edges $e$ of $P$.
Therefore $r(x) = r(x')$, hence $f(x) = f(x')$. This shows that the function $f: S_0 \to \{0,1\}$ is locally constant, hence continuous (in $S_0)$.

\textbf{Step II:} Extend $f$ continuously from $S_0$ to $S_0 \cup S_1 = \bR^2 \setminus P$.

Suppose $p \in S_1$. Let $p_i$ be the unique vertex of $P$ that meets
$R(p)$, i.e., $p \in p^+_i$.
Note that $p \not= p_i$, i.e., $p \in \mbox{ relint } p^+_i$. Let $e_1 = [p_{i-1},p_i], e_2 = [p_i,p_{i+1}]$ be the two edges of $P$ incident with
$p_i$. Define $L = p+\bR v$. $L$ is the vertical line through $p$. Denote by $L^-, L^+$ the two closed half-planes of $\bR^2$ bounded by $L$. None of the edges
$e_1,e_2$ is included in $L$, and they may be either in the same half-plane $L^-$ or $L^+$, or in different half-planes. Choose the notation
so that either $(\alpha)$ $e_1 \subset L^-, e_2 \subset L^+$ (Fig. 3) or $(\beta)$ $e_1 \cup e_2 \subset L^+$ (Fig. 4).

\begin{center}
\end{center}

\begin{center}
 Fig. 3: case $\alpha$ \hspace{1cm} ~~~~~~~~~~~~~~~~~~~~~~Fig. 4: case $\beta$
\end{center}

A glance on Figures 3 and 4 shows that for a point $x$ in the vicinity of $p$, but not lying on $L$, the parity of $r(x)$ is the same in either
side of $L$. Hence we can extend the definition of $f$ to $p$ by defining $f(p)$ to be this parity. To make this into a formal argument consider the
closed set $\triangle =_{\rm def} P \cup ({\rm vert} P \setminus \{p_i\})^+$. This set includes the boundary of $e^+$, for every edge $e$ of $P$, except for
$e^+_1$ and $e^+_2$. It also includes the boundaries of $e^+_1$ and $e^+_2$, except for $p^+_i \setminus \{p_i\}$, and it does not contain the point $p$. Put
$\varepsilon =_{\rm def} \mbox{ dist}(p,\triangle) > 0$, and define $U =_{\rm def} \{x \in \bR^2: \|x-p\| < \varepsilon\} = {\rm int} B^2(p,\varepsilon)$.
Note that if $x \in U$, then the closed interval $[p,x]$ misses $\triangle$. Now make the following observations.
\begin{enumerate}
\item[(I)] If $e$ is any edge of $P$, other than $e_1$ and $e_2$, then the interval $[p,x]$ does not meet the boundary of $e^+$, and therefore $p$ and $x$
are either both in $e^+$, or both not in $e^+$.
\item[(II)] If, say, $e_1 \subset L^-$ and $x \in {\rm int} L^-$ then, moving along the interval $[p,x]$ from $p$ to $x$, we start at a point $p \in p^+_i \subset
{\rm bd} e^+_1$, move into int$e^+_1$, and do not hit the boundary of $e^+_1$ again. Therefore $x \in {\rm int} e^+_1$. The same holds with $L^-$ replaced by
$L^+$, and/or $e_1$ replaced by $e_2$. It follows that in case $(\alpha)$: if $x \in U \setminus L$, then $x$ belongs to exactly one of the sets $e^+_1,e^+_2$. And it
follows that in case $(\beta)$: if $x \in U \cap {\rm int} L^-$, then $x$ belongs to none of the sets $e^+_1,e^+_2$; if $x \in U \cap L^+$, then $x$ belongs to both of them.
\item[(III)] If $p_j \in {\rm vert} P \setminus \{p_i\}$, then $p^+_j \subset \triangle$, and therefore $x \notin p^+_j$, who-ever $x \in U$.
\item[(IV)] If $x \in U \setminus L$, then clearly $x \notin p^+_i$. If $x \in U \cap L$, then the interval $[p,x]$ lies on $L$, contains a point $p \in p^+_i \setminus \{p_i\}$
and does not meet $p_i$; therefore $x \in p^+_i \setminus \{p_i\}$ (= relint$p^+_i$). From these observations we infer:
\begin{enumerate}
\item[(A)] $U \setminus L \subset S_0$ and $f$ is constant on $U \setminus L$.
\item[(B)] $U \cap L \subset S_1$.
\end{enumerate}
\end{enumerate}
Now define $f(p)$ to be the constant value that $f$ takes on $U \setminus L$. Clearly, if we apply the same procedure to any point $p' \in U \cap L$, we will end up
with a value $f(p')$ equal to the value $f(p)$ just defined. (Note that any $\varepsilon'$-neighborhood of $p' \,(\varepsilon' > 0)$ contains points of $U \setminus L$.) Thus we
have extended $f$ to a locally constant, hence continuous function $f: \bR^2 \setminus P \to \{0,1\}$.

To complete the proof of statement (E), we define, as indicated after (F) above, the sets ext$P =_{\rm def} f^{-1}(0)$ and int$P =_{\rm def} f^{-1} (1)$. These are clearly
two disjoint open sets in $\bR^2$, whose union is dom$f = \bR^2 \setminus P$. Note that $\bR^2 \setminus {\rm conv}P \subset {\rm ext}P$ and, therefore, int$P \subset {\rm conv} P$.
Thus ext$P$ is unbounded and int$P$ is bounded.

We still have to show that every point of $P$ is a boundary point of both int$P$ and ext$P$ (and therefore int$P \not= \emptyset, {\rm ext} P \not= \emptyset$).
Since the boundaries of int$P$ and of ext$P$ are closed sets, it suffices to show that the common boundary points of int$P$ and ext$P$ are dense in $P$.

For any vertex $p_i \, (1 \le i \le n)$ the intersection of the vertical line $p_i + \bR v$ with an edge $e$ of $P$ is at most a singleton. Thus $e \setminus \cup \{p_i +
\bR v: 1 \le i \le n\}$ is dense in $e$, and $P \setminus \cup \{p_i + \bR v: 1 \le i \le n\}$ is dense in $P$. If $x \in P \setminus \cup \{p_i + \bR v: 1 \le i \le n\}$, then
$x$ belongs to the relative interior of some edge $e$ of $P$. If $\varepsilon > 0$ is sufficiently small, then the points $x + \varepsilon v, x - \varepsilon v$ are both in
$S_0$, the half-line $R(x + \varepsilon v)$ meets $e$, in addition to all edges met by $R(x-\varepsilon v)$. Thus $r(x+ \varepsilon v) = 1 + r (x-\varepsilon v)$, and $f(x+\varepsilon v)
\not= f(x-\varepsilon v)$, i.e., $\{f(x-\varepsilon v), f(x+\varepsilon v)\} = \{0,1\}$. Thus $x$ is a common boundary point of int$P$ and ext$P$. This finishes the proof
of (E).

\section{Proof of (F)}

Put $I_i =_{\rm def} [p_{i-1},p_i], 1 \le i \le n$, the edges of $P$, and for $i = 1,2, \dots, n$ let $u_i$ be a unit vector perpendicular to aff$I_i$. Choose
the orientation of $u_i$ in such a way that for each point $b \in {\rm relint} I_i$ and for all sufficiently small positive value of $\varepsilon, b + \varepsilon u_i \in
{\rm ext}P$ and $b- \varepsilon u_i \in {\rm int}P$. Define $u_{i,i+1} =_{\rm def} u_i + u_{i+1}, \, 1 \le i \le n$ (the indices are taken modulo $n$, i.e., $p_n = p_0,
u_{n+1} = u_1, u_{n,n+1} = u_{n,1} = u_n+u_1$).

\begin{lemma}\label{lem3.1}
If $\varepsilon$ is a sufficiently small positive number, then $p_i + \varepsilon u_{i,i+1} \in {\rm ext}P$, and $p_i - \varepsilon u_{i,i+1} \in
{\rm int}P$ for $1 \le i \le n$.
\end{lemma}

\textbf{Proof:} The edges $I_i,I_{i+1}$ lie in two rays (half-lines) $L_i,L_{i+1}$ bounded by $p_i$, say $L_i = p_i + \bR^+ v_i, L_{i+1} = p_i + \bR^+ v_{i+1}$, where
$v_i, v_{i+1}$ are suitable unit vectors orthogonal to $u_i$, $u_{i+1}$, respectively.

\begin{center}
\end{center}
\begin{center}
~\phantom{000000}(a) ~~\hspace{4.5cm} (b) \hspace{5cm} (c)

Fig. 5
\end{center}

If $\varepsilon$ is a sufficiently small positive number $(0 < \varepsilon < {\rm dist}(p_i,P \setminus ({\rm relint} (I_i \cup I_{i+1}))$, then $B^2 (p_i,\varepsilon)
\setminus P = B^2 (p_i,\varepsilon) \setminus (L_i \cup L_{i+1})$. The union $L_i \cup L_{i+1}$ divides $B^2 (p_i, \varepsilon)$ into two open sectors, $B^2 (p_i, \varepsilon)
\cap {\rm int} P$ and $B^2 (p_i, \varepsilon) \cap {\rm ext} P$. If $L_i,L_{i+1}$ are collinear $(v_{i+1} = -v_i)$, then each one of these two sectors is an open half disc.
In this case $u_i = u_{i+1}$ (Fig. 5(a)), $u_{i,i+1} = 2u_i = 2u_{i+1}$, and the lemma holds trivially. If $u_i,u_{i+1}$ are not collinear, then one of the sectors
is larger than a half disc, and the other is smaller.
In both cases we have
\begin{equation}\label{eq3}
\langle u_i, v_{i+1} \rangle = \langle u_{i+1}, v_i \rangle = \sin \alpha\,,
\end{equation}
where $\alpha$ is the central angle of the sector $B^2 (p_i,\varepsilon) \cap {\rm ext}P$ at $p_i \, (0 \le \alpha \le 360^o)$.

If $\langle u_i, v_{i+1}\rangle < 0$, then $B^2 (p_i, \varepsilon) \cap {\rm ext}P$ is the larger sector (Fig. 5(b)), and if $\langle u_i,v_{i+1}\rangle > 0$, then
$B^2(p_i,\varepsilon) \cap {\rm int}P$ is the larger sector (Fig. 5(c)). Summing up the equalities
\[
\begin{array}{lll}
u_i & = & \langle u_i, u_{i+1}\rangle u_{i+1} + \langle u_i,v_{i+1}\rangle v_{i+1} \,,\\
u_{i+1} & = & \langle u_{i+1}, u_i\rangle u_i + \langle u_{i+1}, v_i\rangle v_i
\end{array}
\]
and using (\ref{eq3}), we find $(1 - \langle u_i, u_{i+1}\rangle) \, (u_i + u_{i+1}) = \sin \alpha\, (v_i + v_{i+1})$.

If $u_i \not= u_{i+1}$, then $1-\langle u_i,u_{i+1}\rangle > 0$, and
\[
u_{i,i+1} = u_i + u_{i+1} = \frac{\sin \alpha}{1-\langle u_i,u_{i+1}\rangle} \cdot (v_i + v_{i+1})\,.
\]
Thus $u_{i,i+1}$ is a positive [resp., negative] multiple of $v_i + v_{i+1}$ when $\sin \alpha > 0$
[resp., $\sin \alpha < 0$]. In both cases, $u_{i,i+1}$ points towards ext$P$, and $-u_{i,i+1}$ towards int$P$. \hfill \rule{2mm}{2mm}

\begin{lemma}\textbf{{\rm (``Push away from \boldmath$P$''\unboldmath)}}\label{lem3.2}
\begin{enumerate}
\item[(a)] Fix $i, \, 1 \le i \le n$, suppose $b \in {\rm relint} I_i$ and $u$ is a vector satisfying $\langle u,u_i \rangle > 0$. Define
$I^0 =_{\rm def} [b,p_i],I^\varepsilon =_{\rm def} [b+ \varepsilon u, p_i + \varepsilon u_{i,i+1}]$ ($u_i, u_{i+1}$ and $u_{i,i+1} = u_i + u_{i+1}$ denote
the same vectors as in the previous lemma). If $\varepsilon$ is a sufficiently small positive number, then $I^\varepsilon \subset {\rm ext}P$ and $I^{-\varepsilon}
\subset {\rm int}P$. (The required smallness of $\varepsilon$ may depend on the choice of the point $b$ and of the vector $u$.)
\item[(b)] Fix $i,  1 \le i \le n$, and define
$J^0 =_{\rm def}  [p_i,p_{i+1}] = I_{i+1}, J^\varepsilon  =_{\rm def}  [p_i + \varepsilon u_{i,i+1}, p_{i+1} + \varepsilon u_{i+1,i+2}]$.
If $\varepsilon$ is a sufficiently small positive number, then $J^\varepsilon \in {\rm ext}P$ and $J^{-\varepsilon} \in {\rm int}P$.
\end{enumerate}
\end{lemma}

\textbf{Proof:}
\begin{enumerate}
\item[(a)] First note that $I^0$ does not meet any edge of $P$ except $I_i$ and $I_{i+1}$. The same holds for $I^\varepsilon$, provided
\[
|\varepsilon| < \min \left(\frac{1}{2}, \frac{1}{\|u\|}\right) \cdot {\rm dist} \left(I^0, P \setminus({\rm relint} (I_i \cup I_{i+1}))\right)\,.
\]
By Lemma \ref{lem3.1}, $p_i + \varepsilon u_{i,i+1} \in {\rm ext}P$ and $p_i - \varepsilon u_{i,i+1} \in {\rm int}P$, provided $\varepsilon$ is
positive and sufficiently small. To complete the proof, it suffices to show that $I^\varepsilon \cap I_i = \emptyset$ and $I^\varepsilon \cap I_{i+1} = \emptyset$
(for sufficiently small $|\varepsilon|, \, \varepsilon \not= 0$).

As for $I_i : \langle u_i, u \rangle > 0$ (given) and $\langle u_i, u_{i,i+1}\rangle = 1 + \langle u_i, u_{i+1}\rangle > 0$. Therefore, for any $\varepsilon \not= 0$
both endpoints of $I^\varepsilon$ lie (strictly) on the same side of the line aff$I_i$, hence $I_i \cap I^\varepsilon = \emptyset$.

As for
$I_{i+1}$: If $I_{i+1}$ and $I_i$ lie on the same line $(u_i = u_{i+1})$, then the previous argument shows that $I_{i+1} \cap I^\varepsilon = \emptyset$ for
all $\varepsilon \not= 0$ as well. If $u_i \not= u_{i+1}$, consider first the case $\langle u_i, v_{i+1}\rangle < 0$. (Fig. 5(b)). For $\varepsilon > 0, I^\varepsilon$
lies in the open half-plane $\{x \in \bR^2 : \langle u_i, x \rangle > \langle u_i,p_i\rangle\}$, whereas $I_{i+1}$ lies in the closed half-plane $\{x \in \bR^2: \langle
u_i, x \rangle \le \langle u_i, p_i \rangle \}$. Therefore $I^\varepsilon \cap I_{i+1} = \emptyset$. For $\varepsilon < 0$,
\[
\langle u_{i+1}, p_i + \varepsilon u_{i,i+1}\rangle =
\langle u_{i+1}, p_i \rangle + \varepsilon (1 + \langle u_i, u_{i+1}\rangle) < \langle u_{i+1},p_i\rangle\,.
\]
On the other hand, $\langle u_{i+1},b\rangle < \langle u_{i+1},p_i\rangle$ (for any point $b \in {\rm relint} I_i$, since $\langle u_{i+1},v_i \rangle < 0$), and
therefore $\langle u_{i+1},b + \varepsilon u\rangle < \langle u_{i+1},p_i\rangle$ for sufficiently small $|\varepsilon|, \varepsilon \not= 0$. Thus both endpoints of
$I^\varepsilon$ lie on the same open side of the line aff$I_{i+1}$, hence $I^\varepsilon \cap I_{i+1} = \emptyset$.

In the case $\langle u_i, v_{i+1}\rangle > 0$ (Fig. 5(c) above), just repeat the previous argument with the roles of $\varepsilon > 0$ and $\varepsilon < 0$
interchanged.
\item[(b)] The proof is similar to that of (a). First, note that $J^0$ does not meet any edge of $P$ except $I_i,I_{i+1}$ and $I_{i+2}$. The same holds for $J^\varepsilon$, provided
\[
|\varepsilon| < \min \left(\frac{1}{2}, \frac{1}{\|u\|}\right) \cdot {\rm dist } \left(J^0, P \setminus {\rm relint} (I_i \cup I_{i+1} \cup I_{i+2})\right)\,.
\]
By Lemma \ref{lem3.1}, $p_i + \varepsilon u_{i,i+1}, p_{i+1} + \varepsilon u_{i+1,i+2} \in {\rm ext}P$ and $p_i - \varepsilon u_{i,i+1}, p_{i+1} - \varepsilon u_{i+1,i+2} \in
{\rm int} P$, provided $\varepsilon$ is positive and sufficiently small. To complete the proof, it suffices to show that $J^\varepsilon \cap I_i = \emptyset,
J^\varepsilon \cap I_{i+1} = \emptyset$ and $J^\varepsilon \cap I_{i+2} = \emptyset$ (for sufficiently small $|\varepsilon|, \varepsilon \not= 0$).

As for $I_{i+1}\!: \langle u_{i+1},u_{i,i+1}\rangle = 1 + \langle u_{i+1},u_i\rangle > 0$ and $\langle u_{i+1},u_{i+1,i+2}\rangle = 1 +
\langle u_{i+1},u_{i+2} \rangle > 0$. Therefore, for any $\varepsilon > 0$, both endpoints of $J^\varepsilon$ lie on the same open side of the line aff$I_{i+1}$,
hence $I_{i+1} \cap J^\varepsilon = \emptyset$.

As for $I_i$: If $I_{i+1}$ and $I_i$ lie in the same line $(u_i = u_{i+1})$, then the previous argument shows that
$I_i \cap J^\varepsilon = \emptyset$ for all $\varepsilon \not= 0$ as well. If $u_i \not= u_{i+1}$, consider first the case $\langle u_i, v_{i+1}\rangle < 0$ (Fig. 5(b)).

For $\varepsilon > 0, J^\varepsilon$ lies in the open half-plane $\{x \in \bR^2: \langle u_{i+1},x\rangle > \langle u_{i+1},p_i\rangle\}$, whereas $I_i$ lies in the closed
half-plane $\{x \in \bR^2: \langle u_{i+1},x\rangle \le \langle u_{i+1},p_i\rangle\}$.
Therefore, $J^\varepsilon \cap I_i = \emptyset$.

For $\varepsilon < 0$, we have $\langle u_i,p_i + \varepsilon u_{i,i+1}\rangle = \langle u_i, p_i \rangle + \varepsilon (1 + \langle u_i,u_{i+1}\rangle) < \langle u_i,p_i\rangle$.

On the other hand, $\langle u_i, p_{i+1}\rangle < \langle u_i,p_i\rangle$ (since $\langle u_i,v_{i+1}\rangle < 0$), and therefore $\langle u_i, p_{i+1} + \varepsilon u_{i+1,i+2} \rangle
< \langle u_i,p_i\rangle$ for sufficiently small $|\varepsilon|$. Thus both endpoints of $J^\varepsilon$ lie on the same open side of the line aff$I_i$, hence
$J^\varepsilon \cap I_i = \emptyset$.

In the case $\langle u_i, v_{i+1}\rangle > 0$ (Fig. 5(c)), just repeat the previous argument with the roles of $\varepsilon > 0$ and $\varepsilon < 0$ interchanged.

As for $I_{i+2}$: Since the roles of $I_i$ and $I_{i+2}$ are interchangeable, the statement proved above for $I_i$ applies to $I_{i+2}$ as well. \hfill \rule{2mm}{2mm}
\end{enumerate}

\begin{definition}\label{def3.1}
Let $p$ be a point in $\bR^2 \setminus P$ (= ${\rm ext}P \cup {\rm int}P$), and $I$ be an edge of $P$. We say that $p$ \emph{sees} $I$ if, for some point $a \in {\rm relint}\ I, [p,a]\cap
P = \{a\}$.
\end{definition}

\begin{lemma}\label{lem3.3}
Assume $p \in \bR^2 \setminus P$. Then $p$ sees at least one edge of $P$.
\end{lemma}

\textbf{Proof:} Assume, w.l.o.g., that $p \in {\rm ext}P$. Let $q$ be a point in int$P$. Let $U$ be a neighborhood of $q$ that lies entirely in int$P$. Choose a point $q' \in U$
such that the line aff$(p,q')$ does not meet any vertex of $P$. (This condition can be met by avoiding a finite number of lines through $p$.) Then the line segment $[p,q']$ must meet
$P$. Let $a$ be the first point of $P$ on $[p,q']$ (starting from $p$). Then $a$ is a relative interior point of some edge $I$ of $P_i$, and $[p,a]\cap P = \{a\}$. \rule{2mm}{2mm}

\begin{definition}{(poldiam(\boldmath$\cdot$\unboldmath))}:\label{def3.2}
For a set $S \subset \bR^2$ and points $a,b \in S$, denote by $\pi_S (a,b)$ the smallest number of edges of a polygonal path that connects $a$ to $b$ within $S$ ($\pi_S(a,b) =_{\rm def}
\infty$ if no such polygonal path exists). If $S$ is polygonally connected, then $\pi_S (\cdot,\cdot)$ is an integer valued metric on $S$. The \emph{polygonal diameter} of $S$ is defined
as poldiam$(S)=_{\rm def}$ ${\rm sup}\{\pi_S(a,b) : a,b \in S\}$.
\end{definition}

To prove (F) in Section 1 above, it suffices to show that poldiam(int$P$)$< \infty$ and poldiam(ext$P$)$< \infty$.
The following theorem does it.

\begin{theorem}{\rm \textbf{(straightforward upper bound on poldiam(int\boldmath$P$\unboldmath) and poldiam(ext\boldmath$P$\unboldmath))}}\label{theo3.1}
If $P$ is a simple closed $n$-gon $(n \ge 3)$ in $\bR^2$, then we have that {\rm poldiam(int$P$)} and {\rm poldiam(ext$P$)} are both $\le \lfloor\frac{n}{2}\rceil + 3$.
\end{theorem}

\textbf{Proof:} Assume that $a,b$ are two points in the same component (int$P$ or ext$P$) of $\bR^2 \setminus P$. By Lemma \ref{lem3.2}, $a\,[b]$ sees at least one edge
$I'\,[I'']$ of $P$ via $\bR^2 \setminus P$ (possibly $I' = I''$). The set $P \setminus ({\rm relint} (I' \cup I''))$ consists of at most two simple polygonal paths
$P',P''$, the shorter one of which, say $P'$, concatenated by $I',I''$ in both of its endpoints is of the form $\langle J_0,J_1,\dots,J_m,J_{m+1}\rangle$, where $m \le
\lfloor\frac{n-2}{2}\rfloor = \lfloor\frac{n}{2}\rfloor-1, J_0, J_1\dots,J_{m+1}$ are edges of $P \, (\{J_0,J_{m+1}\} = \{I',I''\})$, $J_{i-1}$ and $J_i$ share a vertex $q_i$
for $i=1, \, 2,\dots,m+1$, $a$ sees via $\bR^2 \setminus P$ a point $a' \in {\rm relint} J_0$, and $b$ sees via $\bR^2 \setminus P$ a point $b' \in {\rm relint} J_{m+1}$.

Thus $\langle a,a',q_1,q_2,\dots,q_m,q_{m+1},b',b\rangle$ is a polygonal path of $m + 4 \le \lfloor\frac{n}{2}\rfloor - 1 + 4 = \lfloor\frac{n}{2}\rfloor + 3$ edges that connects
$a$ to $b$ and runs along $P$ except for $[a,a']$ and $[b',b]$. By Lemma \ref{lem3.2}, this path can be pushed away from $P$ into $\bR^2 \setminus P$, thus producing a polygonal
path of $m + 4 \le \lfloor\frac{n}{2}\rfloor +3$ edges that connects $a$ to $b$ via $\bR^2 \setminus P$. \hfill \rule{2mm}{2mm}

\section{Tight upper bounds on poldiam(int\boldmath$P$\unboldmath) and on poldiam(ext\boldmath$P$)}

Theorem \ref{theo3.1} gives a upper bound on poldiam(int$P$)~ [poldiam(ext$P$)] which is somewhat ``naive'', but sufficient to prove (F) in Section 1 above. Here we ``squeeze''
the proof of Theorem \ref{theo3.1} to obtain a tight result.

\begin{theorem}{\rm (Main Theorem)}\label{theo4.1}
Let $P$ be a simple closed $n$-gon in $\bR^2, n \ge 3$. Then
\begin{enumerate}
\item[(a)] the polygonal diameter of {\rm int}$P$ is $\le \lfloor\frac{n}{2}\rfloor$, and the polygonal diameter of {\rm ext}$P$ is $\le \lceil\frac{n}{2}\rceil$;
\item[(b)] for every $n \ge 3$, there is an $n$-gon $P_n$ for which \emph{both} bounds are attained.
\end{enumerate}
\end{theorem}

\textbf{Proof of Theorem 4.1(a):} First note that if $P$ is a convex polygon, then poldiam(int$P)= 1 \le \lfloor\frac{n}{2}\rfloor$, and it can be easily checked that
poldiam(ext$P)=2 \le \lceil\frac{n}{2}\rceil$. (If we consider the closures, however, we find that poldiam(cl int$P)=1$, whereas poldiam(cl ext$P)= 3$ if $P$ has
parallel edges, and equals $2$ otherwise.) This settles the case $n=3$ ($P_3$ is just a triangle). If $n=4$ and $P$ is not convex, then ext$P$ is the union of three convex sets
(two open half-planes and a wedge), each two having a point in common, and therefore poldiam (ext$P) = 2 = \lceil\frac{n}{2}\rceil$. This settles the case $n=4$ for ext$P$.

In view of the proof of Theorem \ref{theo3.1} and the foregoing discussion, we can establish the bounds on poldiam(int$P$) and poldiam(ext$P$) as claimed in Theorem \ref{theo4.1}(a) by showing the following:

\begin{theorem}\label{theo4.2}
Let $P$ be a closed simple $n$-gon in $\bR^2$.
\begin{enumerate}
\item[(i)]
If $n \ge 4$ and $a,b \in {\rm int}P$, then there are two vertices $a',b'$ of $P$ such that a sees $a'$ via int$P$, $b$ sees $b'$ via {\rm int}$P$, and $a',b'$ are at
most $\lfloor\frac{n}{2}\rfloor-2$ edges apart on $P$. (Recall that ``$a$ sees $a'$ via {\rm int}$P$'' means just: $]a,a'[\subset {\rm int}P$.)
\item[(ii)] If $n \ge 5$ and $a,b \in {\rm ext}P$, then there are two vertices $a',b'$ of $P$ such that a sees $a'$ via {\rm ext}$P$, $b$ sees $b'$ via {\rm ext}$P$, and $a',b'$ are at
most $\lceil\frac{n}{2}\rceil-2$ edges apart on $P$,\\
\emph{or:} $\pi_{{\rm ext P}}(a,b) \le 3 \left(\le \lceil\frac{n}{2}\rceil\right.$ for $\left.n \ge 5\right)$.
\end{enumerate}
\end{theorem}

\begin{remark}\label{rem1}
The condition $n \ge 5$ in the first part of Theorem \ref{theo4.2} (ii) cannot be relaxed to $n \ge 4$: Let $P_4 = \langle p_0,p_1,p_2,p_3\rangle$ be a convex
quadrilateral, and let $a,b \in {\rm ext}P_4$,
$a$ close to $[p_0,p_1]$ and $b$ close to $[p_2,p_3]$. Then $a$ and $b$ do not see a common vertex of $P_4$.
\end{remark}

\begin{lemma}\label{lem4.1}
Let $P$ be a simple closed polygon in $\bR^2$. Let $\lceil b',p\rceil$ be an edge of $P$, $a,b$ two points such that $a \in \bR^2 \setminus P$, $b \in ]b',p]$ $($=$[b',p]
\setminus \{b'\})$ and a sees $b$ $($via $\bR^2 \setminus P$$)$. Then a sees $($via $\bR^2 \setminus P$$)$ a vertex of $P$ included in $[a,b',b] \setminus [a,b]$.
\end{lemma}

\textbf{Proof:} If a sees $b'$ then we are done. Otherwise the polygon $P \setminus ]b',p[$ meets the set $[a,b,b'] \setminus [b',b]$. For $0 \le \lambda \le 1$, define
$b(\lambda) =_{\rm def} (1-\lambda) b + \lambda b'$, and let $\lambda_0$ be the smallest value of $\lambda$, $0 \le \lambda \le 1$, such that
$[a,b (\lambda)] \cap (P \setminus ]b',p[) \not= \emptyset$
$(0 < \lambda_0 \le 1; \lambda_0 = 1$ is possible). Let $c'$ be the point of $[a,b(\lambda_0) ] \cap P$ nearest to $a$. Then $c'$ is a vertex of $P$, $c' \in [a,b,b'] \setminus
[a,b]$ and $a$ sees $c'$. \hfill \rule{2mm}{2mm}

\begin{corollary}\label{cor4.1}
Let $P$ be a simple closed $n$-gon, $n \ge 3$, in $\bR^2$. Every point $a \in \bR^2 \setminus P$ sees via $\bR^2 \setminus P$ at least two vertices of $P$.
\end{corollary}

\textbf{Proof:} Let $R$ be a ray emanating from $a$ that meets $P$. By a slight rotation of $R$ around $a$ we may assume that $R$ does not meet any vertex of $P$, but still
$R \cap P \not= \emptyset$. Let $b$ be the first point of $R$ that belongs to $P$ (starting from $a$). By assumption $b \in [b',b''[$ for some edge $[b',b'']$ of $P$. By Lemma \ref{lem4.1}, $a$ sees
via $\bR^2 \setminus P$ a vertex $c'$ $[c'']$ of $P$ included in $[a,b,b'] \setminus [a,b]$ [included in $[a,b,b''] \setminus [a,b]$], and clearly $c' \not= c''$.
\hfill \rule{2mm}{2mm}

\begin{lemma}\label{lem4.2}
Let $P$ be a simple closed $n$-gon, $n \ge 4$, in $\bR^2$, and let $a \in \bR^2 \setminus P$. If every ray emanating from a meets $P$, then a sees via $\bR^2 \setminus P$ two
\emph{non-adjacent} vertices of $P$.
\end{lemma}

\begin{remark}\label{rem2}
The condition that every ray emanating from $a$ meets $P$ is met by every point $a \in {\rm int}P$.
\end{remark}

\textbf{Proof:} By Corollary \ref{cor4.1}, $a$ sees a vertex $c$ of $P$ via $\bR^2 \setminus P$. Consider the ray $R =_{\rm def} \{a + \lambda (a-c) : \lambda \ge 0\}$ that
emanates from $a$ in a direction \emph{opposite} to $c$. By our assumption, $R$ meets $P$. Let $b$ be the first point of $R$ that belongs to $P$. If $b$ is a vertex
of $P$, then a sees the two vertices $b,c$ via
$\bR^2 \setminus P$. These vertices are \emph{not adjacent}, since $[c,b] \cap P = \{c,b\}$. Otherwise, if $b$ is not a vertex of $P$, then $b$ is a relative interior point
of an edge $[b',b'']$ of $P$ $(R \cap ]b',b''[ = \{b\}$). By Lemma \ref{lem4.1}, $a$ sees via $\bR^2 \setminus P$ a vertex $c'$ $[c'']$ of $P$ included in $[a,b,b'] \setminus [a,b]$
[included in $[a,b,b''] \setminus [a,b]$]. Clearly, $c' \not= c''$ and $c',c''$ are non-adjacent in $P$ unless $c' = b'$ and $c'' = b''$. In this
case $a$ sees via $\bR^2 \setminus P$ both couples of vertices $\{c,b'\}$ and $\{c,b''\}$. At least one of these couples is \emph{non-adjacent} in $P$, otherwise $P$ would
be a triangle, contrary to the assumption that $n \ge 4$. \hfill \rule{2mm}{2mm}

\textbf{Proof of Theorem 4.2:}\label{theo4.2}
\begin{enumerate}
\item[(i)] Suppose $P$ is a simple closed $n$-gon, $n \ge 4$, in $\bR^2$. Define $S =_{\rm def} {\rm int}P$, and assume $a,b \in S$. If $n = 4,5$, then cl$S$ (=$P \cup {\rm int}P$) is
starshaped with respect to a vertex of $P$. (If $n = 5$, then $S$ can be triangulated by two interior diagonals with a common vertex.) In this case $a$ and $b$ see via $S$ a
common vertex $a'$ of $P$. Define $b' =_{\rm def} a'$; we find that $a',b'$ are at zero edges apart on $P$. But $0 \le 0 = \lfloor\frac{n}{2}\rfloor-2$ for $n=4,5$.

Assume, therefore, that $n \ge 6$, and that $a$ and $b$ do not see a common vertex of $P$ via $S$. By Lemma \ref{lem4.2}, $a$ sees via $S$ two non-adjacent vertices $a',a''$ of
$P$. These vertices divide $P$ into two paths $P_1,P_2$, each having $\le n-2$ edges. Applying Lemma \ref{lem4.2} again, we find that $b$ sees via $S$ two non-adjacent
vertices $b',b''$ of $P$ and $\{a',a''\} \cap \{b',b''\} = \emptyset$.

If both $b'$ and $b''$ are interior vertices of the same path, say $P_1$, then they divide $P_1$ into three parts. The middle part has at least two edges, and the two
extreme parts together have at most $n-4$ edges. The shorter extreme part, with endpoints (say) $a',b'$, has at most $\lfloor\frac{n-4}{2}\rfloor = \lfloor\frac{n}{2}\rfloor-2$
edges.

If, however, $b'$ is an interior vertex of $P_1$ and $b''$ is an interior vertex of $P_2$, then they divide $P_1$ and $P_2$ into four polygonal paths, each
one of which having one endpoint $b'$ or $b''$. The shortest of these paths has at most $\lfloor\frac{n}{4}\rfloor$ edges. But $\lfloor\frac{n}{4}\rfloor \le \lfloor\frac{n}{2}\rfloor-2$
for $n \ge 6$.
\item[(ii)] Assume $n \ge 5$, define $T = {\rm ext}P$, and let $a,b \in T$. Then either
\begin{enumerate}
\item[(A1)] every ray emanating from $a$ meets $P$, \emph{or}
\item[(A2)] some ray emanating from $a$ misses $P$.
\end{enumerate}
Similarly, either
\begin{enumerate}
\item[(B1)] every ray emanating from $b$ meets $P$, \emph{or}
\item[(B2)] some ray emanating from $b$ misses $P$.
\end{enumerate}
If (A1) and (B1) hold, then both $a$ and $b$ see via $T$ two non-adjacent vertices of $P$ (Lemma \ref{lem4.2}).
If $n \ge 6$, this implies that $a[b]$ sees a vertex $a'$ $[b']$ of $P$ such that $a',b'$ are at most $\lfloor\frac{n-4}{2}\rfloor = \lfloor\frac{n}{2}\rfloor-2 \le
\lceil\frac{n}{2}\rceil-2$ or $\lfloor\frac{n}{4}\rfloor \le \lfloor\frac{n}{2}\rfloor-2 \le \lceil\frac{n}{2}\rceil-2$ edges apart on $P$, as in the proof of part (i) above.
If $n=5$, then $a$ sees via $T$ a vertex $a'$ of $P$, and $b$ sees via $T$ a vertex $b'$ of $P$, where $a'$ and $b'$ are either equal or adjacent, i.e., $a',b'$ are at most one
edge apart on $P$. But for $n=5$ one has $1 \le \lceil\frac{n}{2}\rceil-2$.

If (A2) and (B2) hold, then, due to the compactness of $P$, we can find rays $R_a = \{a + \lambda u : \lambda \ge 0\}$ and $R_b = \{b + \lambda v : \lambda \ge 0\}$ that
miss $P$, where the direction vectors $u$ and $v$ are \emph{linearly independent}. When $\lambda$ is sufficiently large, the segment $[a+\lambda u, b + \lambda u]$ misses $P$.
Therefore $\pi_T (a,b) \le 3 \left(\le \lceil\frac{n}{2}\rceil\right.$ for $\left.n \ge 5\right)$ if $R_a \cap R_b = \emptyset$, and $\pi_T(a,b) = 2 < 3 \left(\le\lceil\frac{n}{2}\rceil \mbox{ for } n \ge 5 \right)$ if $R_a \cap R_b \not= \emptyset$.

If (A1) and (B2) hold, then $a$ sees via $T$ two non-adjacent vertices $a',a''$ of $P$, which divide $P$ into two paths $P_1,P_2$ (with disjoint relative interiors) each one of which
having $\le n-2$ edges. The point $b$,
however, sees two distinct vertices $b',b''$ of $P$, which may be adjacent (Corollary 4.1).
If $\{a',a''\} \cap \{b',b''\} \not= \emptyset$, then again $\pi_T(a,b) \le 2 < 3 \left(\le\lceil\frac{n}{2}\rceil\right.$ for $\left.n \ge 5\right)$. If $\{a',a''\}
\cap \{b',b''\} = \emptyset$, then $b'$ and $b''$ are interior vertices of $P_1$ or $P_2$, or both. If $b'$ and $b''$ belong to different paths, then (as in the proof of part
(i) above) they divide $P_1$ and $P_2$ into four polygonal paths, each having one endpoint $b'$ or $b''$. The shortest one of these paths has at most $\lfloor\frac{n}{4}\rfloor$ edges. But $\lfloor\frac{n}{4}\rfloor \le
\lceil\frac{n}{2}\rceil-2$ for $n \ge 5$. If both $b'$ and $b''$ are interior vertices of the same path, say $P_1$, then (as in the proof of part (i) above) they divide $P_1$
into three parts. The two extreme parts together have at most $n-2-1=n-3$ edges. The shortest extreme part with endpoints (say) $a',b'$ has at most $\lfloor\frac{n-3}{2}\rfloor$
edges. But
$\lfloor\frac{n-3}{2}\rfloor = \lfloor\frac{n-1}{2}\rfloor-1 = \lceil\frac{n}{2}\rceil-2$ for all $n \in \bN $.

The same applies when (A2) and (B1) hold. This finishes the proof of Theorem 4.2. \hfill \rule{2mm}{2mm}

\end{enumerate}

By this also the proof of Theorem 4.1(a) is finished. \hfill \rule{2mm}{2mm}

\textbf{Proof of Theorem 4.1(b):}

We split our examples into two cases, namely even $n$ and odd $n$, $n \ge 3$.

\begin{example}\label{ex4.1} \textbf{\boldmath$n=2m$\unboldmath~ (even), \boldmath$m \ge 2$\unboldmath}.
Figure 6 shows the example for the case $m = 3\, (n=6)$.

\begin{center}
\end{center}

\begin{center}
 Fig. 6: $m = 3 \, (n=6)$
\end{center}

Here we have
$\pi_{{\rm int} P} (a,b) =  m \,(=3) = \lfloor\frac{n}{2}\rfloor$ and
$\pi_{{\rm ext} P} (c,d) =  m \,(=3) = \lceil\frac{n}{2}\rceil $.
One can extend the figure inward beyond vertex $\# 4$.
\end{example}

\begin{example}\textbf{\boldmath$n=2m+1$\unboldmath~ (odd), \boldmath$m \ge 1$\unboldmath.}\label{ex4.2}
Figure 7 shows the example for the case $m = 3$ $(n=7)$

\begin{center}
\end{center}

\begin{center}
Fig. 7: $m = 3 \, (n=7)$
\end{center}

We have
$\pi_{{\rm int} P} (a,b) =  m\, (=3) = \lfloor\frac{n}{2}\rfloor$ and
$\pi_{{\rm ext} P} (c,d) =  m+1\, (=4) = \lceil\frac{n}{2}\rceil$.
Again, one can extend the figure inward beyond vertex $\# 4$.

\end{example}


\begin{thebibliography}{9999999}
%
\bibitem{Aleksandrov}{\sc Aleksandrov, P. S.:} \emph{Combinatorial Topology} (in three volumes), translated from the Russian, 1947.
Kombinatornaya Topologiya, by Harace Komm, Graylock Press, Rochester (1956). Reproduced by Dover Publications,  New York (1998).
%
\bibitem{Bensen}{\sc Benson, R. V.:} \emph{Euclidean Geometry and Convexity}, McGraw-Hill, New York (1966).
%
\bibitem{Bertoglio}{\sc Bertoglio, N., Chuaqui, R.:} An elementary geometric nonstandard proof of the Jordan curve theorem, \emph{Geom. Dedicata}
\textbf{51} (1994), 15-27.
%
\bibitem{Courant}{Courant, R., Robbins, H.:} \emph{What is Mathematics}? 4th ed., Oxford University Press, London (2003).
%
\bibitem{Dostal}{\sc Dost\'{a}l, M., Tindell, R.:} The Jordan curve theorem revisited, \emph{Jahresber. Deutsch. Math.-Ver.} \textbf{80} (1978), 111-128.
%
\bibitem{Hille}{\sc Hille, E.:} \emph{Analytic Function Theory}, Vol. I, Ginn and Company, Boston (1959), (Chelsea, 1973).
%
\bibitem{Kuratowski}{\sc Kuratowski, K.:} \emph{Introduction to Set Theory and Topology}, 2nd ed., Pergamon Press and Polish Scientific Publications, Warsaw (1972).
%
\bibitem{Lawson}{\sc Lawson, T.:} \emph{Topology: A Geometric Approach}, Oxford Graduate Texts in Mathematics, Vol. 9, Oxford University Press, London (2003).
%
\bibitem{Moise}{\sc Moise, E. E.:} \emph{Geometric Topology in Dimension 2 and 3}, Springer, New York (1977).
%
\bibitem{Perles}{\sc Perles, M. A., Martini, H., Kupitz, Y. S.:} A Jordan-Brouwer separation theorem for polyhedral pseudomanifolds, \emph{Discrete Comput. Geom.} \textbf{42}
(2009), 277-304.
%
\bibitem{Thomassen}{\sc Thomassen, C.:} The Jordan-Sch\"onflies theorem and the classification of surfaces, \emph{Amer. Math. Monthly} \textbf{99} (1992),
116-130.
%
\end{thebibliography}
\end{document}